\documentclass[a4paper,12pt]{amsart}
\usepackage{amsmath,amsthm,amssymb,latexsym}
\newtheorem{thm}{Theorem}[section]
\newtheorem{lemma}[thm]{Lemma}

\newtheorem{prop}[thm]{Proposition}
\numberwithin{equation}{section}

\def\a{\alpha}

\def\e{\varepsilon}
\def\H{\mathcal{H}}
\def\L{\mathcal{L}}
\def\o{\omega}
\def\d{\delta}

\def\y{\eta}
\def\k{\kappa}
\def\t{\tau}

\def\O{\Omega}
\def\G{\Gamma}

\def\R{\mathbb R}

\def\Z{\mathbb Z}
\def\00{\infty}
\def\->{\rightarrow}

\def\Ti{T_2^{-1}(D_{\varepsilon_i}\cap Q_{\varepsilon_i})}

\input{amssym.def}

\title[stable phase interfaces in the van~der~Waals--Cahn--Hilliard theory]{stable phase interfaces in the van~der~Waals--Cahn--Hilliard theory}

\author[Y. Tonegawa]{Yoshihiro Tonegawa}

\address[Y. Tonegawa]{Department of Mathematics, Hokkaido University, Sapporo 060-0810 Japan.}
\email{tonegawa@math.sci.hokudai.ac.jp}

\author[N. Wickramasekera]{Neshan Wickramasekera}

\address[N. Wickramasekera]{Department of Pure Mathematics and Mathematical Statistics, University of Cambridge, Cambridge, CB3 0WB, United Kingdom.}
\email{N.Wickramasekera@dpmms.cam.ac.uk}

\date{}


\thanks{Y.T. is partially supported by JSPS Grant-in-aid for scientific research 21340033. Both authors wish to thank Mathematisches Forschungsinstitut Oberwolfach for providing the opportunity to initiate this joint research.}

\begin{document}
\setlength\parskip{5pt}
\begin{abstract}
We prove that any limit-interface corresponding to a  locally uniformly bounded, locally energy-bounded 
sequence of  stable critical points of the van der Waals--Cahn--Hilliard energy functionals with perturbation parameter $\to 0^{+}$  
is supported by an embedded smooth stable minimal hypersurface in low dimensions and an embedded smooth stable minimal hypersurface away from a closed singular set of 
co-dimension $\geq 7$ in general dimensions.  

This result was previously known in case the critical points are local minimizers of energy, in which case the limit-interface is locally area minimizing and its (normalized) multiplicity is 1 a.e. 

Our theorem uses earlier work of the first author establishing stability of the limit-interface as an integral varifold, and relies on a recent general theorem of the second author for its regularity conclusions in the presence of higher multiplicity.
\end{abstract}

\maketitle
\makeatletter
\@addtoreset{equation}{section}%
\renewcommand{\theequation}{\thesection.\@arabic\c@equation}
\makeatother

\section{Introduction}
Let $\Omega \subset {\mathbb R}^{n}$ ($n\geq 2$) be a bounded domain and consider the family of energy functionals $E_{\e}$, $\e \in (0, 1)$, arising in  the van der Waals--Cahn--Hilliard theory of phase transitions (\cite{CahnHilliard}; see also \cite{ModicaMortola}), given by 
\begin{equation}
E_{\e}(u) =\int_{\O}\frac{\e|\nabla u|^2}{2}+\frac{W(u)}{\e}\, d{x},
\label{energy}
\end{equation}
where $u:\O\->\R$ 
belongs to the Sobolev space $H^1(\O)=\{u\in L^2(\O)\,:\,
\nabla u\in L^2(\O)\}$ and $W:\R\->\R^+\cup\{0\}$ is a given $C^{3}$ double-well
potential function with (precisely two) strict minima at $\pm 1$ with $W(\pm 1)=0$. When $\e \to 0^{+}$ with $E_{\e}(u_{\e})$ remaining bounded independently of $\e$, it is clear (from the bound on the second term of the integral above)
that $u_{\e}$ must stay close to $\pm 1$ on a bulk region in 
$\O$ and typically (i.e.\ in case the sets $\{u_{\e}  \approx 1\}$ and $\{u_{\e} \approx  -1\}$ each has measure $\geq$ a fixed proportion of the measure of $\O$ ) there is a transition layer of thickness $O(\e),$
which we may call an ``interface region'' or a ``diffused interface''. 

In the past few decades it has been established that in the presence of a uniform bound on the energy $E_{\e}(u_{\e})$ and under natural variational hypotheses on $u_{\e}$ of varying degrees of generality, for small $\e> 0$, the interface region corresponding to $u_{\e}$ is close to a generalized minimal hypersurface $V$ of $\O$ (the ``limit-interface'' as $\e \to 0^{+}$) and that $E_{\e}(u_{\e})$ approximates a fixed multiple of the $(n-1)$-dimensional area of this hypersurface.
L.~Modica (\cite{Modica}) and P.~Sternberg (\cite{Sternberg}) established this, in the framework of 
$\Gamma$-convergence, for absolutely energy minimizing $u_{\e}$ satisfying a uniform volume constraint; they proved that in this case, the limit-interface $V$ is area minimizing in an appropriate class. R.~Kohn and P.~Sternberg (\cite{KohnSternberg}) studied the locally energy minimizing case, again in the context of $\Gamma$-convergence. More recently, J. Hutchinson and the first author (\cite{HutchinsonTonegawa}) showed that $V$ is a stationary integral varifold  if $u_{\e}$ are assumed to be merely volume-unconstrained critical points of $E_{\e}$ (Theorem~\ref{thm1} below), and that $V$ is an integral varifold with constant generalized mean curvature when the $u_{\e}$ are critical points subject to a volume constraint (see also \cite[Theorem~7.1]{RoegerTonegawa}). Subsequently, the first author (\cite{Tonegawa}) showed that whenever the $u_{\e}$ are unconstrained stable critical points of $E_{\e},$ the limit stationary integral varifold $V$ is stable in the sense that $V$ admits a generalized second fundamental form which satisfies the stability inequality (Theorem~\ref{stabilitythm} below).  

With regard to smoothness of $V$ in the absence of an energy minimizing hypothesis, little has been known beyond the following theorem of the first author (\cite{Tonegawa}): \emph{Suppose that $n=2$, $\e_{i} \to 0^{+}$ as $i \to \infty$ and that for each $i=1, 2, 3, \ldots,$ 
$u_{\e_i}  \in H^{1}(\Omega)$ is a {\rm stable}
critical point of $E_{\e_i}$ with $\sup_{\O} |u_{\e_{i}}| + E_{\e_i}(u_{\e_i})\leq c$ for  some $c>0$ independent of $i$. Then there exists a locally finite union $L$ of non-intersecting lines of $\Omega$ such that after passing to  a subsequence of $\{\e_{i}\}$ without changing notation, for any $0<s<1,$ the sequence of sets $\{{x}\in \O\,:\,
|u_{\e_i}({x})|\leq s \}$ converges  locally in Hausdorff distance to $L$.} Thus in case $n=2,$ any stable diffused interface must be close to 
non-intersecting lines for sufficiently small positive values of the parameter $\e$.

It has remained an open question whether one can make analogous conclusions in dimensions $n>2$. Here we give an affirmative answer to this question in all dimensions. Specifically, we prove (in Theorem \ref{thm2} below) that {\em if $u_{\e_{i}}$ are uniformly bounded stable critical points of $E_{\e_{i}}$ with no volume constraint and with uniformly bounded energy, then for $2 \leq n \leq 7$, there exists an {\em embedded smooth stable minimal hypersurface} $M$ of $\Omega$ such that after passing to a subsequence  of $\{\e_{i}\}$ 
without changing notation, for each fixed $s \in (0, 1)$, the sequence of interface regions $\{{x} \in \Omega \, : \, |u_{\e_{i}}({x})| < s\}$ converges locally in Hausdorff distance to $M$; for $n\geq 8,$ the limit stable minimal hypersurface $M$ may carry an interior singular set,  which is discrete if $n=8$ and has Hausdorff dimension at most $n-8$ if $n \geq 9.$} This regularity result was known
for the limit-interfaces corresponding to sequences $\{u_{\e_{i}}\}$ of energy minimizers since in that case the limit-interfaces are area-minimizing and the well known regularity theory for locally area minimizing currents is applicable.  The new result in this paper is that the stability hypothesis, which is much  
weaker than any energy minimizing hypothesis, suffices to guarantee the same regularity for the limit-interfaces.

The main reason why, in \cite{Tonegawa},  the interface regularity was established only in case $n=2$ and not for $n>2$ was that while in case $n=2$ (i.e.\ when the interface is a 1-dimensional varifold), the structure theorem (due to  W.~Allard and F.~Almgren \cite{AllardAlmgren}) for stationary 1-dimensional varifolds is applicable to the limit-interface, 
there was no sufficiently general regularity theory available at the time 
for higher dimensional stable integral varifolds. In contrast to  limit-interfaces corresponding to sequences of locally energy minimizing critical points of $E_{\e}$, a general stable limit-interface may develop higher multiplicity, {\em a priori} variable even locally. This fact gives rise to significant difficulties that need to be overcome in understanding smoothness properties of a stable limit-interface, and is the reason why the regularity question for stable limit-interfaces in arbitrary dimension remained unresolved prior to the present work. Note that the Schoen--Simon regularity theory (\cite{SchoenSimon}),  which was the most general theory available for stable hypersurfaces at the time when work in \cite{Tonegawa} was carried out, requires knowing {\em a priori} that the singular set (in particular the set of those singular points where the varifold has tangent hyperplanes of multiplicity $\geq 2$) is sufficiently small,  a hypothesis which appears to be difficult to verify directly for a stable limit-interface. The key new input to this problem is the recent work of the second author 
(\cite{Wickramasekera}), which gives a necessary and sufficient geometric structural condition for a general stable
codimension 1 integral varifold to be regular (Theorem~\ref{wick} below). Here we show that the limit-interface in question satisfies precisely this  structural condition; its regularity then follows directly from the general theory of \cite{Wickramasekera}.

While the present work as well as the series of works mentioned above (\cite{Modica,
Sternberg,HutchinsonTonegawa,RoegerTonegawa,Tonegawa}) investigate the 
general character of limit-interfaces, there
have been a number of articles which address the question of existence 
of critical points of \eqref{energy} whose interface regions converge
to a given minimal hypersurface. In this direction we 
mention the work by F.~Pacard and R.~Ritor\'{e} (\cite{PacardRitore}), M.~Kowalczyk 
(\cite{Kowalczyk}) and a number of recent joint works by M.~del~Pino, 
M.~Kowalczyk, F.~Pacard, J.~Wei and J.~Yang (see for example \cite{DPKW}). We refer the reader to the recent survey 
paper by Pacard \cite{Pacard} for a complete list of references.

\section{Hypotheses and the main results}

In this section, we state the hypotheses on $W$ and $u_{\e}$, state our main theorem (Theorem~\ref{thm2}) and recall some definitions and known results needed for its proof. We will give the proof of Theorem~\ref{thm2} in Sections~\ref{sec3} and~\ref{sec4}. 

We assume:

\begin{itemize}
\item[(A1)] $W\in C^3({\mathbb R})$, $W \geq 0$ and $W$ has
precisely three critical points,  two of which are minima at $\pm 1$ with $W(\pm 1)=0$ and $W''(\pm 1)>0$, and the third a local maximum between $\pm 1$.\\
\item[(A2)] $\e_{1}, \e_{2}, \e_{3}, \ldots$ are positive numbers with $\lim_{i\->\00}\e_i=0,$ the constants $c_{1}, c_{2}$ are positive and for each $i=1, 2, 3, \ldots,$ the function $u_{\e_{i}} \in H^{1}(\O)$ and satisfies $E_{\e_i}(u_{\e_i})\leq c_1$ and $\sup_{\O}|u_{\e_i}|\leq c_2$.\\
\item[(A3)] $u_{\e_{i}}$ is a stable critical point of $E_{\e_{i}}$ for each $i=1, 2, 3, \ldots$ 
Thus $u_{\e_{i}}$ solves, weakly,  
\begin{equation}
-\e_i\Delta u_{\e_i}+\frac{W'(u_{\e_i})}{\e_i}=0\hspace{.5cm}\mbox{on $\O$} \label{critical}
\end{equation}
and satisfies
\begin{equation}
 \int_{\Omega} \e_{i}|\nabla \phi|^{2} + \frac{W^{\prime\prime}(u_{\e_{i}})}{\e_{i}}\phi^{2} \geq 0 \hspace{.5cm} \mbox{for each} \;\; \phi\in C^1_c(\O).
\label{secondvar}
\end{equation}
\end{itemize}

\noindent
{\bf Remarks:} {\bf(1)} Hypotheses (\ref{critical}) and (\ref{secondvar}) are equivalent,  respectively, to the conditions 
$$\left.\frac{d}{dt}\right|_{t=0} \, E_{\e_{i}}(u_{\e_{i}} + t\phi) = 0 \;\; {\rm and}$$ 
$$\left.\frac{d^2}{dt^2}\right|_{t=0} E_{\e_i}(u_{\e_i}+t\phi)\geq 0$$ 
for each $\phi \in C^{1}_{c}(\Omega).$ 

\noindent
{\bf (2)} Since $u_{\e_{i}}$ is bounded (by hypothesis 
(A2)), it follows from standard elliptic theory that $u_{\e_{i}} \in C^{3}(\Omega)$ with (\ref{critical})  satisfied pointwise in $\O.$

We shall use the following notation throughout the paper: 

A general point in $\R^n$ is denoted by ${x}=(x_1,\cdots,x_n)$ or by $({y},{z})$ with ${ y}=(y_1,y_2)\in \R^2$
and ${z}=(z_1,\cdots,z_{n-2})\in \R^{n-2}$. For ${x}\in {\mathbb R}^{n}$ and $r>0$, 
we let $B_r({x})=\{\tilde{x}\in \R^n\,:\, |\tilde{x}-{x}|<r\}$ and 
abbreviate $B_r({0})$ as $B_{r}$. For $k \neq n$, ${w}\in \R^k$ and 
$r>0$, we let $B^k_r({w})=\{\tilde{w}\in\R^k\,:\, |\tilde{w}-{w}|<r\}$
and abbreviate $B^{k}_{r}(0)$ as $B_r^k$. 
The $k$-dimensional Lebesgue measure on $\R^k$ will be denoted by ${\mathcal L}^{k}$ and 
$\o_k = \L^k(B_1^k)$. The $k$-dimensional
Hausdorff measure on $\R^n$ is denoted by ${\mathcal H}^{k}$.

In order to characterize the limit-interfaces, we use the notion of varifolds 
(\cite{Allard,Simon}). Let $G_{n-1}(\O)$ be the product space  $\O\times G(n,n-1)$ with
the product topology, where $G(n,n-1)$ is the space
of $(n-1)$-dimensional subspaces in $\R^n$. We identify each element $S \in G(n,n-1)$ with
the $n\times n$ matrix corresponding to the orthogonal projection of ${\mathbb R}^{n}$ onto $S$.
A Radon measure $V$ on $G_{n-1}(\O)$ is called an {\it $(n-1)$-varifold} (henceforth
just called a {\it varifold}). For a varifold $V$, $\|V\|$ shall denote the associated 
{\em weight measure} on $\O$, defined by
\begin{equation*}
\|V\|(\phi)=\int_{G_{n-1}(\O)}\phi({x})\, dV({x,S})\hspace{.5cm}
\mbox{for $\phi\in C_c(\O)$.}
\end{equation*}
The support of the measure $\|V\|$ in $\O$ is denoted by ${\rm spt}\,\|V\|$.
We say that $V$ is {\it integral} if there exists a countably $(n-1)$-rectifiable subset $M$ of $\O$ 
and an $\H^{n-1}$ measurable, positive integer valued function $\theta \, : \, M \to \Z^{+}$ such that $V$ is given by
\begin{equation}
V(\phi)=\int_M\phi({x},T_{x}M)\theta({x})\, d\H^{n-1}({x})
\hspace{1cm}\mbox{for $\phi\in G_{n-1}(\O),$} \label{integral}
\end{equation}
where $T_{x}M$ is the 
approximate tangent space of $M$ at ${x}.$ The function $\theta$ is called the {\it multiplicity} of $V$.

For an $(n-1)$-dimensional $C^{1}$ submanifold $M$ in $\O$, $|M|$ denotes the
multiplicity 1 integral varifold associated with $M,$ as in \eqref{integral} with $\theta \equiv 1$.

We say that $V$ is {\it stationary} if $V$ has zero first variation with respect to area under deformation by any $C^{1}$ vector field of $\Omega$ with compact support, namely (see \cite{Simon}), if

\begin{equation*}
\int_{G_{n-1}(\O)}{\rm tr}\, (\nabla{g}({x})\cdot{S})\, dV({x,S})=0 
\hspace{1cm}\mbox{for all ${g}\in C^1_c(\O;\R^n)$}.
\end{equation*} 

Here $\cdot$ is the usual matrix multiplication and ${\rm tr}$ is 
the trace operator.

For $V\in G_{n-1}(\O)$, let ${\rm reg}\, V\subset\O$ be the set of regular
points of ${\rm spt}\, \|V\| \cap \O$. Thus, ${x}\in {\rm reg}\, V$ if and only if 
$x \in {\rm spt} \, \|V\|$ and there exists some open ball $B_r({x})\subset \O$ such that ${\rm spt}\, \|V\|\cap 
\overline{B_r({x})}$
is a compact, connected, embedded smooth hypersurface with boundary contained in $\partial B_r({x})$. The interior singular set ${\rm sing} \, V$  of $V$ is defined by 
$${\rm sing}\, V=({\rm spt}\, \|V\|\setminus {\rm reg}\, V) \cap \Omega.$$ 
By definition, ${\rm sing}\, V$ is closed in $\Omega$.

To each $u_{\e_i}$ satisfying (A1)-(A3), we associate the varifold $V_{\e_i}\in G_{n-1}(\O)$ defined by
\begin{equation*}
V_{\e_i}(\phi) =\frac{1}{\sigma}\int_{\{|\nabla u_{\e_i}|>0\}}
\phi({x},{\bf I}-{\bf n}_{\e_i}({x})\otimes{\bf n}_{\e_i}({x}))
\frac{\e_i}{2}|\nabla u_{\e_i}|^2\, d{x} \hspace{1in} (\star)
\end{equation*}
for $\phi\in C_c(G_{n-1}(\O))$, where ${\bf n}_{\e_i}({x})=
\frac{\nabla u_{\e_i}({x})}{|\nabla u_{\e_i}({x})|}$, 
$\sigma=\int_{-1}^1\sqrt{W(s)/2}\, ds$, ${\bf I}$ is the $n\times n$ identity
matrix and $\otimes$ is the tensor product. Note that $\|V_{\e_i}\|$
then corresponds simply to $\frac{1}{\sigma}\frac{\e_i}{2}|\nabla u_{\e_i}|^2\, d{x}$.

As a consequence of hypothesis (A2), there exists a subsequence of $\{V_{\e_i}\}_{i=1}^{\00}$ converging as varifolds on $\O$ (i.e.\ as Radon measures on $G_{n-1}(\O)$) to some 
$V \in G_{n-1}(\O)$. Our main result, which concerns regularity of $V$,  is the following:

\begin{thm}
Let the hypotheses be as in (A1)-(A3) and let $V_{\e_{i}}$ be defined by ($\star$). Let 
$V \in G_{n-1}(\O)$ be such that $\lim_{i\to \infty} \, V_{\e_{i}} =V,$ where the convergence is as varifolds on $\O.$ Then ${\rm sing}\, V$
is empty if $2\leq n\leq 7$, ${\rm sing}\, V$ is a discrete set of points if $n=8$
and $\H^{n-8+\gamma}({\rm sing}\, V)=0$ for each $\gamma>0$ if $n\geq 9$; furthermore, ${\rm reg}\,V
=({\rm spt}\, \|V\|\setminus{\rm sing}\, V) \cap \O$ is an embedded smooth stable minimal hypersurface of $\Omega$. 
\label{thm2}
\end{thm}

\noindent
{\bf Remarks:} {\bf (1)} As just mentioned, if the hypotheses are as in (A1)-(A3), then after passing to a subsequence of $\{\e_{i}\}$ without changing notation, we obtain $V \in G_{n-1}(\O)$ such that $V_{\e_{i}} \to V$ as varifolds on $\O.$ It is of course possible that $V = 0$. 

\noindent
{\bf (2)} There exists $u_{0} \in BV(\O)$ such that after passing to a subsequence of $\{\e_{i}\}$ without changing notation, 
$u_{\e_{i}} \to u_{0}$ in $L^{1}(\O);$ in fact $u_{0}(x) = \pm 1$ for a.e. $x \in \O,$  and hence the sets 
$\{u_{0} = 1\}$ and $\{u_{0} = -1\}$ have finite perimeter in $\O$ (see the discussion in \cite{HutchinsonTonegawa}, pp. 51-52) and 
${\rm spt} \,\|\partial \, \{u_{0} = 1\}\| \cap \O \subset {\rm spt} \, \|V\|$ (see \cite{HutchinsonTonegawa}, Theorem 1). Thus in particular, 
$V \neq 0$ is implied by the condition that $u_{0}  \not\equiv 1$ on $\O$ and $u_{0} \not\equiv -1$ on $\O.$
  
We now further discuss known results we shall need for the proof of Theorem~\ref{thm2}. 

The following theorem, due to Hutchinson and the first author (\cite{HutchinsonTonegawa}), says among other things that a limit varifold $V$ corresponding to a sequence of critical points of $E_{\e}$ 
(i.e.\ $V \in G_{n-1}(\O)$ arising as the varifold limit of a sequence $\{V_{\e_{i}}\},$ where $\e_{i} \to 0^{+}$ and $V_{\e_{i}} \in G_{n-1}(\O)$ is defined by ($\star$) with $u_{\e_{i}}$ satisfying \eqref{critical}) is a stationary integral varifold, and that ${\rm spt} \, \|V\|$ indeed is the limit-interface corresponding to $\{u_{\e_{i}}\}$ in the sense made precise in part (3) of the theorem. Note that no stability hypothesis is necessary for this result.

\begin{thm}{\rm (}\cite[Theorem 1 and Proposition 4.3]{HutchinsonTonegawa}{\rm )}
Suppose that (A1), (A2) and \eqref{critical} hold, and let $V \in G_{n-1}(\O)$ be such that $V = \lim_{i\to \infty} \, V_{\e_{i}},$ where $V_{\e_{i}}$ is as in ($\star$) and the convergence is as varifolds on $\O.$ Let $U$ be an open subset of $\O$ such that the closure of $U$ is contained in $\O.$ Then
\begin{itemize}
\item[(1)] $V$ is a stationary integral varifold on $\O$.
\item[(2)] 
\begin{equation*}
\lim_{i\->\00}\int_U\left|\e_i\frac{|\nabla u_{\e_i}|^2}{2}-
\frac{W(u_{\e_i})}{\e_i}\right|\, d{x}=0.
\label{discrepancy}
\end{equation*}
\item[(3)] For each $s \in (0, 1)$, $\{|u_{\e_i}|\leq s\}\cap U$
converges to ${\rm spt}\, \|V\|\cap U$ in Hausdorff distance.
\end{itemize}
\label{thm1}
\end{thm}

In order to discuss the additional known results relevant to us concerning \emph{stable} critical points of 
$E_{\e}$, it is convenient to introduce the following notation: For $u\in C^2(\O),$ let ${\mathbf B}_{u}$ be the non-negative function defined by  
\begin{equation}
{\mathbf B}_{u}^2 =\frac{1}{|\nabla u|^2}\sum_{i,j=1}^n
u_{x_i x_j}^2-\frac{1}{|\nabla u|^4}\sum_{i=1}^n\left(\sum_{j=1}^n
u_{x_j}u_{x_i x_j}\right)^2
\label{defsec}
\end{equation}
on $\{|\nabla u|>0\}$ and ${\mathbf B}_{u} =0$ on $\{|\nabla u|=0\}$. Here and subsequently, $u_{x_{i}},$  $u_{x_{i}x_{j}}$ denote the partial derivatives $\frac{\partial \, u}{\partial \, x_{i}},$ $\frac{\partial^{2} \, u}{\partial x_{j} \partial x_{i}}$ respectively. Note that the expression on the right 
hand side of (\ref{defsec}) is non-negative when $\nabla \, u \neq 0$ and is invariant under orthogonal transformations of ${\mathbb R}^{n}$.

We have the following:

\begin{lemma}{\rm (}\cite{PadillaTonegawa}{\rm )}\label{stab}
Let $\e \in (0, 1)$, $u \in C^{2}(\O)$ and suppose that $u$ is a stable critical point of $E_{\e}$ in the sense that \eqref{critical} and \eqref{secondvar} are satisfied with $\e$ in place of $\e_{i}$ and $u$ in place of $u_{\e_{i}}.$ Then 
\begin{equation}
\int_{\O}{\mathbf B}_{u}^2 |\nabla u|^2\phi^2\, d{x}
\leq \int_{\O} |\nabla\phi|^2|\nabla u|^2\, d{x}
\label{secondest1}
\end{equation}
for each $\phi\in C^1_c(\O).$ 
\end{lemma}

One proves \eqref{secondest1} by taking $|\nabla u|\phi$ in place of $\phi$
in the inequality \eqref{secondvar} and  utilizing equation \eqref{critical}. See~\cite{PadillaTonegawa} or ~\cite{Tonegawa} for  details. 

Let the hypotheses be as in (A1)-(A3), and write 
$${\mathbf B}_{\e_i} = {\mathbf B}_{u_{\e_{i}}}.$$
In view of hypothesis (A2), Lemma~\ref{stab} implies that the $L^1$-norm of $\e_i{\mathbf B}_{\e_i}^2|\nabla u_{\e_i}|^2$ is locally uniformly bounded.
Let $\nu$ be a subsequential limit (as Radon measures on $\O$) of the sequence $\e_{i}{\mathbf B}_{\e_{i}}^{2}\left|\nabla u_{\e_{i}}\right|^{2} dx.$ Thus after re-indexing,

\begin{equation}
\nu(\phi)=\lim_{i\->\00}\int_{\O}\e_i\phi{\mathbf B}_{\e_i}^2|\nabla u_{\e_i}|^2\, d{x}
\label{secondest2}
\end{equation}  

for $\phi\in C_c(\O)$.\\

The following crucial stability inequality is established in \cite{Tonegawa}:

\begin{thm}{\rm (}\cite[Theorem 3]{Tonegawa} {\rm )} \label{stable}
Let the hypotheses be as in (A1)-(A3) and let $V_{\e_{i}}$ be defined by ($\star$). Let 
$V \in G_{n-1}(\O)$ be such that $V = \lim_{i\to \infty} \, V_{\e_{i}},$ where the convergence is as varifolds on $\O.$ Then $V$ has a generalized
second fundamental form ${A}$ with its length $|{A}|$ satisfying 
\begin{equation}
\int_{\O}|{A}|^2\phi^2\, d\|V\|\leq \int_{\O}|\nabla \phi|^2\, d\|V\|
\label{stability}
\end{equation}
for all $\phi\in C^1_c(\O)$.
\label{stabilitythm}
\end{thm}

We refer the reader to \cite[Section 2]{Tonegawa} for the definition of the generalized second
fundamental form of a varifold. See also  \cite{Hutchinson}
where the notion was defined originally. 

We end this section with the following elementary consequence of \eqref{critical} which we shall need later:

\begin{lemma} If $u \in C^{2}(\O)$ and \eqref{critical} holds with $\e \in (0, 1)$ in place of $\e_{i}$ and $u$ in place of $u_{\e_{i}}$, then 
\begin{equation}
\left|\nabla\left(\e\frac{\left|\nabla u\right|^2}{2}-\frac{W(u)}{\e}\right)\right|
\leq \e\sqrt{n-1}\,|\nabla u|^2{\mathbf B}_{u}\hspace{.2in} \mbox{in $\O$}.
\label{discrepancyest}
\end{equation}
\end{lemma}

\noindent
{\it Proof}. For any $n\times n$ symmetric matrix $M$ and any unit vector
${\bf m} \in {\mathbb R}^{n}$, one has that  
\begin{equation*}
|M\cdot {\bf m}-({\rm tr}\,
M){\bf m}|\leq \sqrt{n-1}({\rm tr}\, (M^2)-
{\bf m}^t\cdot M^2\cdot{\bf m})^{\frac12}.
\end{equation*}
Using this and \eqref{critical}, we see that on the set $\{\nabla u \neq 0\}$, 
\begin{equation*}
\begin{split}
&\left|\nabla\left(\e\frac{\left|\nabla u\right|^2}{2}-\frac{W(u)}{\e}\right)\right|
=\e\left|\nabla^2u\cdot \nabla u-\Delta u\nabla u\right|\\
&\leq \e\sqrt{n-1}\, |\nabla u|\left({\rm tr}\, ((\nabla^2 u)^2)-\frac{1}{|\nabla u|^2}
\nabla u^t\cdot (\nabla^2 u)^2\cdot\nabla u\right)^{\frac12}=\e\sqrt{n-1}\,|\nabla u|^2
{\mathbf B}_{u}.
\end{split}
\end{equation*}
If $\nabla u = 0$, the inequality holds trivially.
\hfill{$\Box$}\\

\section{Regularity of stable codimension 1 integral varifolds and the proof of the main theorem}\label{sec3}

In this section we recall the main content (Theorem~\ref{wick} below) of the regularity theory of the second author (\cite{Wickramasekera}) for stable codimension 1 integral varifolds and show how it implies our main result (Theorem~\ref{thm2}) concerning regularity of the limit-interfaces, modulo verification of a certain structural condition satisfied by the limit-interfaces. This structural condition is precisely given in Proposition~\ref{mainprop} below, and is necessary to apply Theorem~\ref{wick}. We shall establish its validity in the next section.

Fix an integer $m \geq 2$ and $\a\in (0,1)$. Denote by ${\mathcal S}_{\a}$ the collection of all integral 
$m$-varifolds $V$ on the open unit ball $B^{m+1}_{1} \subset {\mathbb R}^{m+1}$ with ${0}\in {\rm spt}\, \|V\|$, $\|V\|(B^{m+1}_{1})<\00$ and satisfying the following conditions:
\begin{itemize}
\item[(${\mathcal S}\,1$)] (Stationarity) $V$ is a critical point of the $m$-dimensional area 
functional in $B^{m+1}_{1}$, viz.\ $V$ is a stationary integral $m$-varifold on $B^{m+1}_1$.\\
\item[(${\mathcal S} \,2$)] (Stability) $V$ satisfies
\begin{equation}
\int_{{\rm reg}\, V}|{A}|^2\phi^2\, d\H^{m}\leq
\int_{{\rm reg}\,V}|\nabla^{{\rm reg} \, V}\phi|^2\, d\H^{m}
\label{stabreg}
\end{equation}
for all $\phi\in C^1_c({\rm reg}\, V)$, where ${A}$
denotes the (classical) second fundamental form of ${\rm reg}\, V,$ $|{A}|$ its length and 
$\nabla^{{\rm reg}\, V}$ is the gradient operator on 
${\rm reg} \, V$.\\
\item[(${\mathcal S} \, 3$)] ($\a$-Structural Hypothesis) For each ${x}\in {\rm sing}\, V$,
there exists no $\rho>0$ such that ${\rm spt}\, \|V\|\cap B_{\rho}({x})$
is equal to the union of a finite number of $m$-dimensional embedded $C^{1,\a}$
submanifolds-with-boundary of $B_{\rho}({x})$ all having common boundary
in $B_{\rho}({x})$ equal to an $(m-1)$-dimensional embedded $C^{1,\a}$ submanifold
of $B_{\rho}({x})$ containing ${x}$, and no two intersecting except
along their common boundary.
\end{itemize}

With these hypotheses, we have the following:

\begin{thm}{\rm (}\cite[Theorem 3.1]{Wickramasekera}{\rm )}
If $V\in {\mathcal S}_{\a},$ then ${\rm sing}\, V\cap B^{m+1}_1$ is empty if $2\leq m \leq 6$,
${\rm sing} \, V \cap B_{1}$ is a discrete set of points if $m=7$  and $\H^{m-7+\gamma}({\rm sing}\, V\cap B_1)=0$ for each $\gamma>0$ if $m\geq 8$.
\label{wick}
\end{thm}

An obvious yet extremely useful feature of Theorem~\ref{wick} is that it suffices, when applying the theorem, to verify the 
$\alpha$-Structural Hypothesis for points $x \in {\rm sing} \, V$ in the complement of a set $Z$ having 
$(m-1)$-dimensional Hausdorff measure zero; we do not need to know that such $Z$ is closed. We rely on this fact in an essential way in the present application, in which the $\alpha$-Structural Hypothesis is verified by way of the following proposition. We shall  prove this proposition in the next section.

\begin{prop} Let $V$ be as in Theorem~\ref{thm2}. There exists a (possibly empty) Borel
set $Z\subset{\rm spt}\,\|V\| \cap \O$ with $\H^{n-2}(Z)=0$ such that for each ${x}\in ({\rm spt}\, \|V\|\setminus Z) \cap \O$ and each tangent cone ${\mathbf C}$ to $V$
at ${x}$, ${\rm spt}\, \|{\mathbf C}\|$ is not equal to a union of three or 
more half-hyperplanes of $\R^n$ meeting along an $(n-2)$-dimensional
affine subspace. 
\label{mainprop}
\end{prop}

Theorem~\ref{thm2} follows directly from Proposition~\ref{mainprop} and Theorem~\ref{wick}.

\noindent
\emph{Proof of Theorem~\ref{thm2}}. 
Let $V$ be as in Theorem~\ref{thm2}, $y \in {\rm spt} \, \|V\| \cap \Omega$ and 
$\rho \in (0, {\rm dist} \, (y, \, \partial \, \O))$. In order to prove Theorem~\ref{thm2}, it clearly suffices to establish its conclusions with $\widetilde{V} = \eta_{y, \rho \; \#} \, V$ in place of $V$, which we can achieve by Theorem \ref{wick} if we can show that $\widetilde{V} \in {\mathcal S}_{\alpha}$ for some $\alpha \in (0, 1)$. Here and subsequently, $\eta_{y, \rho} \, : \, {\mathbb R}^{n} \to {\mathbb R}^{n}$ is the map defined by $\eta_{y, \rho}(x) = \rho^{-1}(x - y)$ and $\eta_{y, \rho \, \#}$ denotes the push-forward of $V$ by $\eta_{y, \rho}.$ 

It follows from Theorem \ref{thm1} that $\widetilde{V}$ satisfies (${\mathcal S} \, 1$). To verify that $\widetilde{V}$ satisfies (${\mathcal S} \, 2$), note the following two facts: (i) on the regular part of the varifold, the generalized second fundamental form  in \eqref{stability} coincides with the classical second fundamental 
form (\cite{Hutchinson}); (ii) by the constancy theorem (\cite{Simon}, Theorem 41.1), the multiplicity function of $V$  is constant on each connected component of ${\rm reg}\, V$ so we can replace $d\|V\|$ in \eqref{stability}  by $d\H^{n-1}$ whenever $\phi \in C^{1}_{c}(\Omega \setminus {\rm sing} \, V)$. Thus given $\phi \in C^{1}_{c} \, ({\rm reg} \, \widetilde{V})$, 
we may choose any extension $\widetilde{\phi} \in C^{1}_{c}(B_{1} \setminus {\rm sing} \, \widetilde{V})$
of $\phi$ such that $\nabla \widetilde{\phi} = \nabla^{{\rm reg} \, V} \phi$ on ${\rm reg} \, V$ and 
use \eqref{stability} with $\widetilde{\phi} \circ \eta_{y, \, \rho}$ in place of $\phi$ to deduce that 
$\widetilde{V}$ satisfies (${\mathcal S} \, 2$).

To verify that $V$ (and hence also $\widetilde{V}$) satisfies (${\mathcal S} \, 3$) (for any 
$\alpha \in (0, 1)$), 
let $\alpha \in (0, 1)$ and suppose that there is a point ${x}\in {\rm spt}\, 
\|V\|\cap\O$ and $\rho \in (0, {\rm dist}({x},\partial \O))$ such that ${\rm spt}\,
\|V\|\cap B_{\rho}({x})$ is a union of three or more $C^{1,\a}$ 
hypersurfaces-with-boundary meeting along a common $(n-2)$-dimensional $C^{1,\a}$
submanifold $L$ of $B_{\rho}({x})$ with ${x}\in L$. It is standard to
see, with the help of the Hopf boundary point lemma for divergence form
elliptic operators (\cite[Lemma 7]{FinnGilbarg}, see also \cite[Lemma 10.1]{HardtSimon})
that at every point along $L$, these hypersurfaces-with-boundary must meet
transversely. Hence the unique tangent cone to $V$ at any $\tilde{x}\in L$
is supported by a union of three or more half-hyperplanes meeting along a common
$(n-2)$-dimensional subspace. For any $\widetilde{x} \in L\setminus Z$, this directly contradicts Proposition~\ref{mainprop}, where $Z$ is the set as in Proposition~\ref{mainprop}. Note that 
$L \setminus Z \neq \emptyset$ since 
$\H^{n-2}(Z)=0$ by Proposition~\ref{mainprop}. Thus $V$ must satisfy (${\mathcal S} \, 3$).
Hence $\widetilde{V} \in {\mathcal S}_{\alpha}$, and Theorem \ref{thm2} follows from Theorem \ref{wick}.
\hfill{$\Box$}

\section{Structural condition for the stable limit-interfaces}\label{sec4}

To complete the proof of Theorem~\ref{thm2}, it only remains to give a proof of Proposition~\ref{mainprop}, which we shall do in this section.  

Let the hypotheses be as in (A1)-(A3) and let $\nu$ be the Radon measure on $\O$ defined by \eqref{secondest2}. Let $V$ be as in Theorem~\ref{thm2}, obtained possibly after passing to a suitable subsequence of $\{\e_{i}\}$ and the corresponding subsequence of $\{u_{\e_{i}}\}.$ Let 
\begin{equation*}
Z =\left\{{x}\in {\rm spt}\,\|V\| \cap \O\, :\, \limsup_{r\-> 0}\frac{\nu(B_r({x}))}{r^{n-3}}
>0\right\}.
\end{equation*}
It is standard to see that $\H^{n-3+\gamma}(Z)=0$ for each $\gamma>0$ 
and thus in particular that $\H^{n-2}(Z)=0$. We will show that Proposition~\ref{mainprop} holds 
with this $Z.$ 

To obtain a contradiction assume that we have a point ${x}\in{\rm spt}\,
\|V\|\setminus Z$ where a tangent cone ${\mathbf C}$ to $V$ has the property that ${\rm spt}\,
\|{\mathbf C}\|$ is equal to a union of three or more half-hyperplanes meeting along a common 
$(n-2)$-dimensional subspace $S({\mathbf C})$. Without loss of generality we may assume
that ${x}={0}$ and that $S({\mathbf C}) = \{0\} \times {\mathbb R}^{n-2}.$ Thus we may write
\begin{equation*}
{\rm spt}\, \|{\mathbf C}\|= \cup_{j=1}^{N} P_{j}
\end{equation*}
for some $N \geq 3,$ where for each $j=1, 2, \ldots, N,$ 
$$P_{j} =  \{t{\mathbf p}_{j} \, : \, t \geq 0\} \times {\mathbb R}^{n-2}$$ 
with ${\bf p}_1,\cdots,{\bf p}_N\in \R^2$ distinct vectors such that
$|{\bf p}_{j}|=\frac12$ for $j=1, 2, \ldots, N.$  

By the definition of tangent cone, there exists a sequence $r_i\->0$ 
such that $\lim_{i\->\00}\y_{r_i\,\#}V={\mathbf C}$. Here $\eta_{r}$ is the map 
$x \mapsto r^{-1}x.$ Since $V_{\e_i}\-> V$,
we may choose  a subsequence of $\{\e_{i}\}$ for which, after relabeling, we have that  
$\lim_{i\->\00}\y_{r_i\,\#} \, V_{\e_i}={\mathbf C}$ and $\lim_{i\->\00}\frac{\e_i}{r_i}
=0$. Letting $\tilde{\e}_i=\frac{\e_i}{r_i}$ and defining  
$\tilde{u}_{\tilde{\e}_i}(\tilde{x})=u_{\e_i}(r_i\tilde{x})$, we then have that
\begin{equation*}
\y_{r_i\,\#} \, V_{\e_i}(\phi)=\frac{1}{\sigma}\int_{\{|\nabla \tilde{u}_{\tilde{\e}_{i}}| > 0\}} \phi(\tilde{x},{\bf I}
-{\bf n}_{\tilde{\e}_i}\otimes{\bf n}_{\tilde{\e}_i})\frac{\tilde{\e}_i}{2}
|\nabla \tilde{u}_{\tilde{\e}_i}|^2\, d{\tilde{x}}
\end{equation*}
for $\phi\in C_c(G_{n-1}(r_i^{-1}\O))$, where $\tilde{\mathbf n}_{\tilde{\e}_{i}} = \frac{\nabla \tilde{u}_{\tilde{\e}_{i}}}{|\nabla \tilde{u}_{\tilde{\e}_{i}}|}$. Since $\lim_{i\->\00}\frac{\nu(B_{2r_i})}{(2r_i)^{n-3}}=0$, we may choose a further subsequence of $\{\e_{i}\}$ without changing notation such that
\begin{equation*}
\lim_{i\->\00}\frac{1}{(2r_i)^{n-3}}\int_{B_{2r_i}}\e_i{\mathbf B}_{\e_i}^2
|\nabla u_{\e_i}|^2\, d{x}=0.
\label{secdec}
\end{equation*}
With the change of variables as above, this is equivalent to 
\begin{equation*}
\lim_{i\->\00}\int_{B_2}\tilde{\e}_i{\tilde{\mathbf B}}_{\tilde{\e}_i}^2|\nabla
\tilde{u}_{\tilde{\e}_i}|^2\, d\tilde{x}=0,
\end{equation*}
where ${\tilde{\mathbf B}}_{\tilde{\e}_i}$ is defined by \eqref{defsec} with $\tilde{u}_{\tilde{\e}_i}$
in place of $u$. 

By \cite[Prop. 3.4]{HutchinsonTonegawa}, for each open set $U$ with 
$U \subset\subset \O$, there exist constants $c = c(c_{1}, n, {\rm dist}\, (U, \partial \, \O))$ and 
$\e_{0} = \e_{0}(c_{1}, n, {\rm dist} \, (U, \partial \, \O))$ such that if $B_{r}({x_{0}}) \subset U$ and 
$s \in (0, r]$, then for all $i$ sufficiently large to ensure $\e_{i} \leq \e_{0}$, 
\begin{eqnarray}\label{monotonicity}
&&r^{1-n} \int_{B_{r}(x_{0})} \left(\e_{i}\frac{|\nabla u_{\e_{i}}|^{2}}{2} + \frac{W(u_{\e_{i}})}{\e_{i}}\right) \, dx - s^{1-n} \int_{B_{s}(x_{0})} \left(\e_{i}\frac{|\nabla u_{\e_{i}}|^{2}}{2} + \frac{W(u_{\e_{i}})}{\e_{i}}\right) \, dx\nonumber\\ 
&&\hspace{.5in} \geq \int_{s}^{r} \left(\t^{-n}\int_{B_{\t}(x_{0})} \left(\frac{W(u_{\e_{i}})}{\e_{i}} - \frac{\e_{i}}{2}|\nabla u_{\e_{i}}|^{2}\right)^{+} \, dx\right) \, d\t -cr\nonumber\\ 
&&\hspace{2.5in} + \, \e_{i}\int_{B_{r}(x_{0}) \setminus B_{s}(x_{0})} \frac{\left((y - x_{0}) \cdot \nabla u_{\e_{i}}\right)^{2}}{|y - x_{0}|^{n+1}} \, dy.
\end{eqnarray}
This implies in particular that 
$$\int_{B_{2}} \frac{\tilde{\e}_{i}}{2}|\nabla \tilde{u}_{\tilde{\e}_{i}}|^{2} + \frac{W(\tilde{u}_{\tilde{\e}_{i}})}{\tilde{\e}_{i}}  \leq C$$ 
for all sufficiently large $i$, where $C$ is a positive constant depending only on $n$, $c_{1}$ and 
$c_{2}.$ 

Thus hypotheses (A1)-(A3) are satisfied with $\tilde{\e}_{i}$ in place of $\e_{i}$, 
$\tilde{u}_{\tilde{\e}_{i}}$ in place of $u_{{\e}_{i}}$ and $C$ in place of $c_{1},$ 
so by replacing the original sequences $\{\e_{i}\},$ $\{u_{\e_{i}}\}$ with the new sequences 
$\{\tilde{\e}_{i}\},$ $\{\tilde{u}_{\tilde{\e}_{i}}\}$ and the constant $c_{1}$ with $C$, we have 
that (A1)-(A3) hold with $\O = B_{2}$,  together with the additional facts that 
\begin{equation}\label{cone}
\lim_{i\->\00}V_{\e_i}={\mathbf C}
\end{equation}
where the convergence is as varifolds on $B_{2}$ and  that
\begin{equation}\label{secdec2}
\lim_{i\->\00}\int_{B_2}{\e}_i{{\mathbf B}}_{{\e}_i}^2|\nabla
{u}_{{\e}_i}|^2 =0.
\end{equation}

\emph{For the rest of the discussion we shall assume that $W$, $\{\e_{i}\}$, $\{u_{\e_{i}}\}$ satisfy (A1)-(A3) with $\O = B_{2},$ as well as \eqref{cone} and \eqref{secdec2}}. 

Our goal is to derive a contradiction.

\bigskip

\begin{lemma}
Set $c_3=\frac12 \sqrt{\min_{|t|\leq\frac34}W(t)}>0$ and let
\begin{equation*}
D_{\e_i}=\left\{{z}\in B_1^{n-2}\,:\, |\nabla u_{\e_i}({y},{z})|
\geq \frac{c_3}{\e_i}\mbox{ holds for all }{y}\in B_1^2\mbox{ with }|u_{\e_i}
({y},{z})|\leq \frac12\right\}.
\end{equation*}
Then we have that
\begin{equation*}
\lim_{i\->\00}\L^{n-2}(B_1^{n-2}\setminus D_{\e_i})=0.
\end{equation*}
\label{lemma1}
\end{lemma}

\noindent
{\bf Remark:} Note that $D_{\e_i}$ contains the set $D^{\prime}_{\e_{i}}$ of points ${z} \in B_{1}^{n-2}$ where
$|u_{\e_i}({y},{z})|>\frac12$ for all ${y}\in
B_1^2$. We shall prove that $D^{\prime}_{\e_{i}}$ is small in Lemma \ref{lemma2} below.\\

\noindent
{\it Proof}. For each $i=1, 2, 3, \ldots,$ let $\{B_{\e_i}^{n-2}({z}_{i,j})\}_{j=1}^{J_i}$ 
be a maximal pairwise disjoint collection of balls such that 
${z}_{i,j}\in B_1^{n-2}\setminus D_{\e_i}$ for $j=1,\cdots,J_i.$
Then $B_1^{n-2}\setminus D_{\e_i}\subset 
\cup_{j=1}^{J_i}B_{2\e_i}^{n-2}({z}_{i,j})$ and 
by the definition of $D_{\e_i},$ there exists, for each $j=1, 2, \ldots, J_{i}$, a point  ${y}_{i,j}\in B_1^2$
such that $|u_{\e_i}({y}_{i,j},{z}_{i,j})|\leq\frac12$ and
$|\nabla u_{\e_i}({y}_{i,j},{z}_{i,j})|<\frac{c_3}{\e_i}$.
By standard elliptic estimates we have that 
$$\sup_{B_{15/8}} \, |\nabla^{2} \, u_{\e_{i}}| \leq C \sup_{B_{2}} \, \left(|u_{\e_{i}}| + \frac{|W^{\prime}(u_{\e_{i}})|}{\e_{i}^{2}}\right)$$
where $C = C(n)$, whence in view of the hypothesis $\sup_{B_{2}} \, |u_{\e_{i}}| \leq c_{2}$, there exists a fixed number $r_{0} \in (0, 1],$ depending only 
on $n$, $W$ and $c_2,$ such that
\begin{equation}
|u_{\e_i}({x})|\leq \frac34\mbox{ and }|\nabla u_{\e_i}({x})|<
\frac{2c_3}{\e_i}
\label{small}
\end{equation}
for each ${x}\in B_{2r_0 \e_i}({y}_{i,j},{z}_{i,j})$. On this ball we have
\begin{equation}
v_{\e_i}\equiv \frac{W(u_{\e_i})}{\e_i}-\frac{\e_i|\nabla u_{\e_i}|^2}{2} \geq \frac{1}{\e_i}\left(\min_{|t|\leq \frac34}W(t)-2c_3^2\right) \geq \frac{2c_3^2}{\e_i}.
\label{big}
\end{equation}
Since $B^2_{r_0 \e_i}({y}_{i,j})\times\{{z}\}\subset B_{2r_0\e_i}({y}_{i,j}, {z}_{i, j})$ 
for each ${z}\in B_{r_0 \e_i}^{n-2}({z}_{i,j})$, we have by \eqref{big} that
\begin{equation}
2c_3^2 \sqrt{\pi}\, r_0\leq \left(\int_{B_{r_0\e_i}^2 ({y}_{i,j})}
(v_{\e_i}({y},{z}))^2\, d{y}\right)^{\frac12}
\label{big2}
\end{equation}
for each ${z}\in B_{r_0 \e_i}^{n-2}({z}_{i,j}).$ On the other hand, by the relevant 2-dimensional Sobolev inequality and \eqref{discrepancyest},
\begin{equation}
\begin{split}
&\left(\int_{B^2_{r_0\e_i}({y}_{i,j})}(v_{\e_i}({y},{z}))^2\, 
d{y}\right)^{\frac12} \leq \left(\int_{B^2_{1+r_0\e_i}}(v_{\e_i}({y},{z}))^2
\, d{y}\right)^{\frac12}\\
&\hspace{.5in} \leq C\int_{B^{2}_{1+r_{0}\e_{i}}}|v_{\e_{i}}(y, z)| + |\nabla_{y} v_{\e_i}({y},{z})|\, d{y}\\
&\hspace{.5in} \leq C\int_{B^{2}_{1+r_{0}\e_{i}}} |v_{\e_{i}}(y, z)| \, dy + 
C\sqrt{n-1}\e_{i}\int_{B^2_{1+r_{0}\e_{i}}}{\mathbf B}_{\e_i}(y, z)|\nabla u_{\e_i}(y, z)|^2\, d{y}
\end{split}
\label{big3}
\end{equation}
where $C$ is the relevant Sobolev constant. Combining \eqref{big2} and \eqref{big3} we obtain that
\begin{equation*}
c_3^2\sqrt{\pi}r_0\leq  C\int_{B^{2}_{1+r_{0}\e_{i}}} |v_{\e_{i}}(y, z)| \, dy + 
C\sqrt{n-1}\e_{i}\int_{B^2_{1+r_{0}\e_{i}}}{\mathbf B}_{\e_i}(y, z)|\nabla u_{\e_i}(y, z)|^2\, d{y}
\end{equation*}
for all ${z}\in B^{n-2}_{r_0 \e_i}({z}_{i,j})$. Integrating this over
$B_{r_0\e_i}^{n-2}({z}_{i,j})$ first, summing over $j$  and using the fact that 
$\{B_{r_0\e_i}^{n-2}({z}_{i,j})\}_{j=1}^{J_i}$ are pairwise disjoint, we obtain with the help of the Cauchy-Schwarz inequality that 
\begin{equation*}
\begin{split}
&c_3^2\sqrt{\pi}r_0^{n-1}\o_{n-2}\e_i^{n-2}J_i\leq C\int_{B^{2}_{1 + r_{0}\e_{i}} \times B^{n-2}_{1 + r_{0}\e_{i}}} |v_{\e_{i}}| \, dx + C\sqrt{n-1}\e_{i}\int_{B^2_{1 + r_{0}\e_{i}}\times B^{n-2}_{1 + r_{0}\e_{i}}}{\mathbf B}_{\e_i}|\nabla u_{\e_i}|^2\, d{x}\\
&\hspace{.5in} \leq C\int_{B_{15/8}}|v_{\e_{i}}| \, dx + C\sqrt{n-1}\left(\int_{B_{15/8}}\e_{i}|\nabla \, u_{\e_{i}}|^{2} \, dx \right)^{1/2}\left(\int_{B_{15/8}}\e_{i}{\mathbf B}_{\e_{i}}^{2}|\nabla u_{\e_{i}}|^{2} \, dx\right)^{1/2}\\
&\hspace{.5in} \leq C\int_{B_{15/8}}|v_{\e_{i}}| \, dx + Cc_{1}\sqrt{n-1}\left(\int_{B_{15/8}}\e_{i}{\mathbf B}_{\e_{i}}^{2}|\nabla u_{\e_{i}}|^{2} \, dx \right)^{1/2}.
\end{split}
\label{big4}
\end{equation*}
Since $\L^{n-2}(B_1^{n-2}\setminus D_{\e_i})\leq J_i\o_{n-2}(2\e_i)^{n-2}$,
the lemma follows from Theorem~\ref{thm1}(2) and \eqref{secdec2}. 
\hfill{$\Box$}\\

Choose $\d>0$ small enough depending on ${\mathbf C}$ so that 
$\{B_{2\d}^2({\bf p}_i)\}_{i=1}^{N}$ are disjoint, and define

\begin{equation*}
Q_{\e_i}=\left\{{z}\in B_1^{n-2}\,:\, \forall t\in [-1/2,1/2 ],\,
\forall j\in\{1,\cdots,N\},\, \exists {y}\in B_{\d}^2({\bf p}_j)\, 
\mbox{ s.t. }u_{\e_i}({y},{z})=t\right\}.
\end{equation*}

The next lemma is obtained by re-examining \cite[Section 5]{HutchinsonTonegawa}. 

\bigskip

\begin{lemma}
With $Q_{\e_{i}}$ as above, we have that 
\begin{equation*}
\lim_{i\->\00}\L^{n-2}(B_1^{n-2}\setminus Q_{\e_i})=0.
\end{equation*}
\label{lemma2}
\end{lemma}

\noindent
{\it Proof}. For $j=1, 2, \ldots, N$, let 
\begin{equation*}
Q_{\e_i,j} =\left\{{z}\in B^{n-2}_1\,:\, \forall t\in [-1/2,1/2],\,
\exists {y}\in B^2_{\d}({\bf p}_j)\mbox{ s.t. }u_{\e_i}({y},{z})=t\right\}.
\end{equation*}
It suffices to prove that $\lim_{i\->\00}\L^{n-2}(B_1^{n-2}\setminus Q_{\e_i,j})=0$ for each $j=1, 2, \ldots, N$. Thus without loss
of generality, we may assume $j=1$, $P_1=\{y_2=0\}\cap\{y_1 \geq 0\}$ and
${\bf p}_1=(1/2,0)$. With these coordinates, ${\rm spt}\, \|{\mathbf C}\|
\cap B_{\d}({\bf p}_1,{z})=\{y_2=0\}\cap B_{\d}({\bf p}_1,{z})$ for
each ${z}\in \R^{n-2}$. On $B_{\d}({\bf p}_1,{z})$ with ${z}\in
B^{n-2}_1$, $V_{\e_i}$ converge to $\theta_1 |P_1|$ as varifolds, where
$\theta_1 \in {\mathbb N}$ is the multiplicity of ${\mathbf C}$ on $P_1$. By Theorem \ref{thm1}, the
sets $B_{\d}({\bf p}_1,{z})\cap\{|u_{\e_i}|\leq \frac12\}$ converge to 
$P_1 \cap B_{\d}({\bf p}_1,{z})$ in Hausdorff distance. 
Note that $u_{\e_i}({x})$ converges to different values ($\pm 1$)
uniformly on $\left(\cup_{z \in B_{1}^{n-2}} \, B_{\d}({\bf p}_1,{z})\right)\cap\{y_2 >\frac{\d}{2}\}$ and
$\left(\cup_{z \in B_{1}^{n-2}} \, B_{\d}({\bf p}_1,{z})\right) \cap \{y_2<-\frac{\d}{2}\}$ in 
case $\theta_1$ is odd, and to the same value if $\theta_{1}$ is even 
(see the discussion in \cite[p. 78]{HutchinsonTonegawa}). Hence if $\theta_{1}$ is odd, 
by continuity
of $u_{\e_i}$, the function $y_2 \mapsto u_{\e_i}(\frac12,y_2,{z})$ as $y_2$ ranges over $[-\frac{\d}{2},
\frac{\d}{2}]$ takes all values between $-\frac12$ and $\frac12$, so that in 
this case we see that $B^{n-2}_1= Q_{\e_i,1}$ for all sufficiently 
large $i,$ proving the lemma. 

If $\theta_1$ is even, we need to utilize
results in \cite[Section 5]{HutchinsonTonegawa}. 
Assume without loss of generality that
$u_{\e_i}$ converges to $+1$ on both sides of $\{y_2>0\}$ and $\{y_2<0\}$
on $B_{\delta}({\bf p}_1,z)$. Let
\begin{equation*}
\begin{split}
&\hat{B}_{\delta/2}({\bf p}_1,z)=\{(\hat{y}_1,
\hat{y}_2,\hat{z})\,:\, (\hat{y}_1-1/2)^2+|\hat{z}-z|^2<(\delta/2)^2,
\,\, |\hat{y}_2|<\delta/2\} \;\; {\rm and}\\
&S_i=\left\{x\in B_{\delta/2}({\bf p}_1,z)\cap
P_1\,:\, \exists t\in [-1/2, 1/2],\,\, \mbox{with} \,\, \{u_{\e_i}=t\}\cap
T_1^{-1}(x)\cap \hat{B}_{\delta/2}({\bf p}_1,z)=\emptyset\right\}.
\end{split}
\end{equation*}
Here $T_1$ is the orthogonal projection ${\mathbb R}^n\rightarrow \{y_{2} = 0\}$.
By the continuity of the $u_{\e_i}$'s and their local uniform convergence to $+1$ away 
from $P_1$, we have for any $b \in (0, 1/2)$ that $S_i\subset S_i^b\cup \hat{S}_i^b$, where
\begin{equation*}
\begin{split}
S_i^b=&\left\{x\in B_{\delta/2}({\bf p}_1,z)\cap
P_1 \,:\, \{u_{\e_i}\leq 1-b\}
\cap T_1^{-1}(x)\cap \hat{B}_{\delta/2}({\bf p}_1,z)=\emptyset\right\} \;\; {\rm and}\\
\hat{S}_i^b=&\left\{x\in B_{\delta/2}({\bf p}_1,z)\cap
P_1\,:\, \inf_{T_1^{-1}(x)\cap \hat{B}_{\delta/2}({\bf p}_1,z)}
u_{\e_i}\in [-1/2,1-b] \right\}.
\end{split}
\end{equation*}

We claim that for any given sufficiently small $s>0,$ we can choose
small $b = b(s, W)>0$ such that 
\begin{equation}
\limsup_{i\rightarrow\infty}{\mathcal L}^{n-1}(S_i^b)
\leq c(\sigma,n,\theta_1)s.
\label{Si1}
\end{equation}
To see this, we argue as follows: Note first that for any given $s \in (0,1)$ we have the estimates (5.5)-(5.8) of
\cite{HutchinsonTonegawa} with $B_{\delta}
({\bf p}_1,z)$ in place of $B_{3}$, where $b = b(s, W) >0$ in (5.8) is given by 
\cite[Prop. 5.1]{HutchinsonTonegawa}.  
Choose $\eta = \eta(s, W, \delta, \theta_{1}) \in (0, 1)$ and $L = L(s, W) \in (1, \infty)$ as in \cite[Prop. 5.5, 5.6]{HutchinsonTonegawa} with $R = \delta$, $N = \theta_{1},$ and define $G_i$ by
\begin{equation}
\begin{split}
G_i=&\hat{B}_{\delta/2}({\bf p}_1,z)
\cap\{|u_{\e_i}|\leq 1-b\}\cap\\
&\left\{x\,: \, \int_{B_r(x)}\left(\left|\frac{\e_i|\nabla u_{\e_i}|^2}{2}
-\frac{W(u_{\e_i})}{\e_i}\right|+(1-\nu_{2, \, i}^2)\e_i|\nabla u_{\e_i}|^2\right)
\leq \eta r^{n-1}\,\mbox{ if }\, 4\e_i L\leq r\leq \delta\right\}
\end{split}
\label{Gi}
\end{equation}
where $\nu_{2, \, i} = |\nabla u_{\e_{i}}|^{-1}\frac{\partial \, u_{\e_{i}}}{\partial y_{2}}$ if 
$\nabla \, u_{\e_{i}} \neq 0$ and $\nu_{2, \, i} = 0$ otherwise. With the help of the Besicovitch covering theorem and (\ref{monotonicity}) we see that\\
\begin{equation}
\|V_{\e_i}\|(\hat{B}_{\delta/2}({\bf p}_1,z)\cap\{|u_{\e_i}|\leq 1-b\}\setminus
G_i)+{\mathcal L}^{n-1}(T_1(\hat{B}_{\delta/2}({\bf p}_1,z)\cap \{|u_{\e_i}|\leq
1-b\}\setminus G_i))\rightarrow 0.
\label{Gigoto0}
\end{equation}

We also note that there exists $c=c(b,s,W)$ with the following property:
for a.e. $t \in [-1+b, 1 - b],$ the level set $\{u_{\e_{i}} = t\}$ is an 
$(n-1)$-dimensional $C^3$ surface and for any $x\in G_i$ with $u_{\e_i}(x)=t$, the set 
$\{u_{\e_i}=t\} \cap B_{L\e_{i}}(x)$ is a graph $y_2=f(y_1,z_1,\cdots,z_{n-2})$ of a $C^{1}$ function 
$f \, : \, T_{1}(\{u_{\e_{i}} =t\} \cap B_{L\e_{i}}(x)) \to {\mathbb R}$ 
with $|\nabla f| \leq c\eta^{1/(n+1)}$ on $T_{1}(\{u_{\e_{i}} =t\} \cap B_{L\e_{i}}(x))$. This follows from the proof of \cite[Prop. 5.6]{HutchinsonTonegawa}, which yields that for any $x = (y_{1}^{\star}, y_{2}^{\star}, z^{\star})\in G_{i}$, the function
$u_{\e_i}$ in the neighborhood $B_{L\e_{i}}(x)$ is $C^{1}$ close to $\pm q((y_2-y_{2}^{\star})/\e_i),$ where $q$ is the standing wave solution defined by the ODE $q''=W'(q)$ with $q(\pm \infty)=\pm 1$; specifically, letting $\tilde{u}_{\e_i}
(\tilde{y}_1,\tilde{y}_2,\tilde{z})=u_{\e_i}(\e_i\tilde{y}_1+y_1^{\star},
\e_i\tilde{y}_2+y_2^{\star},\e_i\tilde{z}+z^{\star})$ and $\tilde{q}
(\tilde{y}_1,\tilde{y}_2,\tilde{z})=\pm q(\tilde{y}_2+c)$ (so that $q(c)=t$), we have that
$\|\tilde{u}_{\e_i}-\tilde{q}\|_{C^1(B_L)}\leq c\eta^{1/(n+1)}$. 
In particular, we choose $\eta = \eta(s, W, \delta, \theta_{1}) \in (0, 1)$ so small that 
\begin{equation}
\sqrt{1+|\nabla f|^2}\leq 1+s.
\label{gradsmall}
\end{equation}

For $x\in P_1\cap B_{\delta/2}
({\bf p}_1,z)$ and $|t|\leq 1-b$, define $Y^i_x(t)=T_1^{-1}(x)
\cap G_i\cap\{u_{\e_i}=t\}$. We claim that the cardinality $\# Y_x^i(t)$ of $Y_x^i(t)$ is $\leq \theta_1$. To see this, assume for a contradiction that $\# Y_x^i (t)\geq
\theta_1+1$, and let $Y'$ be any subset of $Y_x^i (t)$ such that
$\# Y'=\theta_1+1$. Then
we have by \cite[Prop. 5.6]{HutchinsonTonegawa}
\begin{equation}
I \equiv \sum_{\tilde{x}\in Y'}\frac{1}{\omega_{n-1}(L\e_i)^{n-1}}
\int_{B_{L\e_i}(\tilde{x})}\frac{\e_i|\nabla u_{\e_i}|^2}{2}+\frac{W(u_{\e_i})}{\e_i}
\geq (\# Y')(2\sigma -s)
\label{Iineq1}
\end{equation}
while by \cite[Prop. 5.5]{HutchinsonTonegawa},
\begin{equation}
\begin{split}
I&\leq s+\frac{1+s}{\omega_{n-1}\delta^{n-1}}
\int_{\{\tilde{x}\,|\, {\rm dist}\,(Y',\tilde{x})<\delta\}}
\frac{\e_i|\nabla u_{\e_i}|^2}{2}+\frac{W(u_{\e_i})}{\e_i}\\
&\leq s+\frac{1+s}{\omega_{n-1}\delta^{n-1}}2\|V_{\e_i}\|
(B_{\delta+o(1)}(x))+o(1)
\end{split}
\label{Iineq2}
\end{equation}
where $o(1)\rightarrow 0$ as $i\rightarrow\infty$ uniformly in $x
\in B_{\delta/2}({\bf p}_1,z)\cap P_1$. Since 
\begin{equation}
\|V_{\e_i}\|
(B_{\delta+o(1)}(x))\leq \sigma\theta_1\delta^{n-1}\omega_{n-1}
+o(1),
\label{Iineq3}
\end{equation}
having $\# Y'= \theta_1 +1$ would contradict \eqref{Iineq1}-\eqref{Iineq3}
for all sufficiently large $i,$ provided $s>0$ is smaller than a number
depending only on
$\theta_1,\, n,\,\delta$ and $\sigma$. We may of course assume that $s>0$ is smaller than
this number to begin with. 

Now defining $w^i$ as in 
\cite[page 52]{HutchinsonTonegawa}, we have by \cite[(5.8)]{HutchinsonTonegawa}
and \eqref{Gigoto0} that 
\begin{equation}
\begin{split}
(\delta/2)^{n-1}\omega_{n-1}\theta_1\sigma&=\lim_{i\rightarrow\infty}
\int_{\hat{B}_{\delta/2}({\bf p}_1,z)}|\nabla w^i|\\
&\leq s+\liminf_{i\rightarrow\infty}\int_{G_i}|\nabla u_{\e_i}|\sqrt{W(u_{\e_i})/2}.
\end{split}
\label{inty1}
\end{equation}
Using the co-area formula, \eqref{gradsmall}, the fact that $\# Y^i_x(t)\leq \theta_1$ 
and \eqref{inty1}, we see that 
\begin{equation}
\begin{split}
(\delta/2)^{n-1}\omega_{n-1}\theta_1\sigma&\leq s+\liminf_{i\rightarrow\00}
\int_{-1+b}^{1-b}{\mathcal H}^{n-1}(\{u_{\e_i}=t\}\cap G_i)\sqrt{W(t)/2}\,
dt\\
&\leq s+\sigma\theta_1(1+s)\liminf_{i\rightarrow\00}{\mathcal L}^{n-1}(T_1(G_i)).
\end{split}
\label{t3}
\end{equation}
Note that $T_1(G_i)\cap S^b_i=\emptyset$ by the definition of $G_{i}$ and $S_{i}^{b}$, and hence 
${\mathcal L}^{n-1}(T_1(G_i))\leq \omega_{n-1}(\delta/2)^{n-1}-{\mathcal L}^{n-1}
(S_i^b)$. In view of \eqref{t3}, this implies that $\limsup_{i\rightarrow\00}
{\mathcal L}^{n-1}(S_i^b)\leq c(\sigma,n,\theta_1)s$, completing the proof of
\eqref{Si1}. 

We next verify that 
\begin{equation}
\hat{S}_i^b\subset T_1(\{|u_{\e_i}|\leq 1-b\}\cap\hat{B}_{\delta/2}(
{\bf p}_1,z)\setminus G_i)
\label{Si2}
\end{equation}
as follows: For any $x=(\hat{y}_{1},0,\hat{z})\in \hat{S}_i^b$, there exist $\hat{y}_2$ with $|\hat{y}_2|\leq \delta/4$ and $t\in [-1/2,1-b]$ with 
$u_{\e_i}(\hat{y}_1,\hat{y}_2,\hat{z})=t$. If $(\hat{y}_1,\hat{y}_2,\hat{z})\in G_i$, 
again as above we have by \cite[Prop. 5.6]{HutchinsonTonegawa} that $u_{\e_i}$ is $C^1$ close to 
$q((y_2-\hat{y}_{2})/\e_i)$  in the $L\e_i$-neighborhood of 
$(\hat{y}_1,\hat{y}_2,\hat{z})$. In particular, we would then 
have $T_1^{-1}(x)\cap\{u_{\e_i}=-3/4\}\cap\hat{B}_{\delta/2}({\bf p}_1,z)\neq \emptyset$, contradicting the assumption that $x \in \hat{S}_i^b$. Thus
$(\hat{y}_1,\hat{y}_2,\hat{z})\in \{|u_{\e_i}|\leq 1-b\}\cap\hat{B}_{\delta/2}(
{\bf p}_1,z)\setminus G_i$, proving \eqref{Si2}. 

It follows from  \eqref{Si2} and \eqref{Gigoto0} that 
\begin{equation}
\lim_{i\rightarrow\00}{\mathcal L}^{n-1}(\hat{S}_i^b)=0.
\label{Si3}
\end{equation}
Since $S_i\subset S_i^b\cup\hat{S}_i^b$, it follows from \eqref{Si1}, \eqref{Si3} and
arbitrariness of $s>0$ that  
\begin{equation}
\lim_{i\rightarrow\00}{\mathcal L}^{n-1}(S_i)=0.
\label{Sif}
\end{equation}

Now to complete the proof, assume contrary to the assertion of the lemma that 
\begin{equation*}
\limsup_{i\rightarrow\00}{\mathcal L}^{n-2}(B_1^{n-2}\setminus
Q_{\e_i,1})>0.
\end{equation*}
Then for some $z\in B_1^{n-2},$ we must have
\begin{equation}
\limsup_{i\rightarrow\00}{\mathcal L}^{n-2}(B_{\delta/4}^{n-2}(z)
\setminus Q_{\e_i,1})>0.
\label{Sig}
\end{equation}
Take any $z'\in B_{\delta/4}^{n-2}(z)\setminus Q_{\e_i,1}$. For any
$y_1$ with $|y_1 - 1/2|<\delta/4$, we have
\begin{equation}
x=(y_1,0,z')\in S_i;
\label{Sig2}
\end{equation}
for if not, there would exist $y_1$ with $|y_1 - 1/2|<\delta/4$
such that $x=(y_1,0,z')\notin S_i$ so that $u_{\e_i}(y_1,y_2,z')$ must take
all values $t\in [-1/2,1/2]$ as $y_2$ ranges over $[-\delta/2,\delta/2]$. 
But this would mean that $z'\in Q_{\e_i,1},$ contrary to our assumption. Thus the claim
\eqref{Sig2} holds, and says that 
\begin{equation}
Z_i \equiv \{(y_1,0,z')\,:\, z'\in B_{\delta/4}^{n-2}(z)\setminus Q_{\e_i,1},
\,\, |y_1 - 1/2|<\delta/4\}\subset S_i.
\label{Sig3}
\end{equation}
But then since ${\mathcal L}^{n-1}(Z_i)=\frac{\delta}{2}{\mathcal L}^{n-2}(B_{\delta/4}^{n-2}(z)\setminus Q_{\e_i,1})$, the statements \eqref{Sif}, \eqref{Sig} and \eqref{Sig3} are contradictory, completing the 
proof of the lemma.
\hfill{$\Box$}\\

\bigskip
In Lemma~\ref{lemma3} below we shall make hypotheses and use notation as follows: Let $u\in C^2(B^2_1\times B^{n-2}_1)$ and suppose that $t$ is a regular 
value of $u$ with $M = u^{-1}(t) \neq \emptyset,$ so that $M$ is an $(n-1)$-dimensional embedded $C^2$ submanifold of 
$B^{2}_{1} \times B^{n-2}_{1}$. Let ${\mathbf B}_{u}$ be the function defined by \eqref{defsec}.
Let $L$ be the set of points $z \in B_{1}^{n-2}$ satisfying the following two requirements: (i) $T_{2}^{-1}(z) \cap M \neq \emptyset$ and (ii) $z$ is a regular value of the map $\left.T_{2}\right|_{M} \, : \, M \to {\mathbb R}^{n-2}$, where $T_{2} \, : \, {\mathbb R}^{n} \to {\mathbb R}^{n-2}$ is the orthogonal projection. Then for each $z \in L,$ $\ell_{z} \equiv T_2^{-1}({z})\cap M$ is  a $C^{2}$ 1-manifold. For $z \in L,$ let $\k_{z}(p)$ denote the curvature of $\ell_{z}$ at $p \in \ell_{z}$. 

\bigskip

\begin{lemma}
Let the hypotheses and notation be as described in the preceding paragraph. Then we have for any Borel set $G\subset L,$
\begin{equation*}
\int_G\, d{z}\int_{\ell_{z}}|\k_{z}|\, ds\leq
\left(\int_{M\cap T_2^{-1}(G)}{\mathbf B}_{u}^2 \, d\H^{n-1}\right)^{\frac12}
\left(\H^{n-1}(M\cap T_2^{-1}(G))\right)^{\frac12}.
\end{equation*}
\label{lemma3}
\end{lemma}

{\it Proof}. Since $L$ is open in $B^{n-2}_1$, we may choose an 
increasing sequence of open sets $L_i\subset\subset L$ such that
$\cup_{i=1}^{\00}L_i=L$. Then $M_i=T_2^{-1}(\overline{L_i})\cap M$ is a 
$C^2$ submanifold-with-boundary $\partial M_i=T_2^{-1}(\partial L_i)
\cap M$ in $B_1^2\times B^{n-2}_1$. On $M_i$, $(u_{y_1},u_{y_2})\neq
(0,0)$. Hence for each point $x \in M_{i}$, there exists $\rho_{x} >0$ such that $M \cap B_{\rho_{x}}(x)$ is the graph of a $C^{2}$ function $v$ defined over an open subset $U$ of either the plane $\{y_{1} = 0\}$ or 
the plane $\{y_{2} = 0\}.$ We may cover $M_i$ with a finite number of such 
coordinate charts $M_{i} \cap B_{\rho_{x}}(x)$. Let us now fix such a chart, and assume without loss of generality that $U \subset \{y_{2} = 0\}$ for that chart, so that $v = v(y_{1}, z)$ for $(y_{1}, z) \in U$ and by the definition of $M,$ $v$ satisfies $u(y_{1}, v(y_{1}, z), z) =t$ for all $(y_{1}, z) \in U.$ In particular, for each $z \in T_{2}(U)$, we have that 
$\ell_{z} \cap B_{\rho}(x) = \{(y_{1}, v(y_{1}, z), z) \, : \, y_{1} \in U \cap T_{2}^{-1}(z)\}$ and hence  
$\k_{z}=v_{y_1 y_1}(1+v_{y_1}^2)^{-\frac32}.$ Using the identity $u(y_{1}, v(y_{1}, z), z) \equiv t$ on $U$, this can be expressed in terms of $u$ as 
\begin{equation*}
\k_{z} = -\frac{u_{y_2 y_2}u_{y_1}^2-2u_{y_1 y_2}u_{y_1}u_{y_2}+u_{y_1 y_1}u_{y_2}^2}
{(u_{y_1}^2+u_{y_2}^2)^{\frac32}}.
\end{equation*}
Since the length element $ds$ is given by 
\begin{equation*}
ds = \sqrt{1+v_{y_1}^2}\, dy_1=\frac{\sqrt{u_{y_1}^2+u_{y_2}^2}}{|u_{y_2}|}\, dy_1,
\end{equation*}
we have that
\begin{equation*}
|\k_{z}|ds=\frac{|u_{y_2 y_2}u_{y_1}^2-2u_{y_1 y_2}u_{y_1}u_{y_2}+u_{y_1 y_1}
u_{y_2}^2|}{|u_{y_2}|(u_{y_1}^2+u_{y_2}^2)}\, dy_1.
\end{equation*}
Next for unit vector ${\bf m}\in {\mathbb R}^n$ with ${\bf m}\perp \nabla u$
and $M=(\nabla^2 u)$, we have 
\begin{equation*}
({\bf m}^t M{\bf m})^2\leq |M{\bf m}|^2
={\bf m}^t M^2 {\bf m}\leq {\rm tr}(M^2)-\frac{1}{|\nabla u|^2}
(\nabla u)^t M^2 (\nabla u)=|\nabla u|^2|{\bf B}_u|^2.
\end{equation*}
Note that here we used the non-negativity of the eigenvalues of $M^2$.
Since $(u_{y_2},-u_{y_1},0,\cdots,0)\perp \nabla u$, we deduce that
\begin{equation*}
|u_{y_2 y_2}u_{y_1}^2-2u_{y_1 y_2}u_{y_1}u_{y_2}+u_{y_1 y_1}
u_{y_2}^2|^2\leq (u_{y_1}^2+u_{y_2}^2)^2{\mathbf B}_{u}^2 |\nabla u|^2,
\end{equation*}
which implies that 
\begin{equation*}
\frac{|u_{y_2 y_2}u_{y_1}^2-2u_{y_1 y_2}u_{y_1}u_{y_2}+u_{y_1 y_1}
u_{y_2}^2|}{|u_{y_2}|(u_{y_1}^2+u_{y_2}^2)} \leq \frac{{\mathbf B}_{u}|\nabla u|}{|u_{y_2}|}.
\end{equation*}
We also note that
\begin{equation*}
\sqrt{1+v_{y_1}^2+v_{z_1}^2+\cdots+v_{z_{n-2}}^2}=\frac{|\nabla u|}{|u_{y_2}|},
\end{equation*} 
so that 
\begin{equation*}
\frac{|u_{y_2 y_2}u_{y_1}^2-2u_{y_1 y_2}u_{y_1}u_{y_2}+u_{y_1 y_1}
u_{y_2}^2|}{|u_{y_2}|(u_{y_1}^2+u_{y_2}^2)} \leq {\mathbf B}_u 
\sqrt{1+v_{y_1}^2+v_{z_1}^2+\cdots+v_{z_{n-2}}^2},
\end{equation*}
where both the expression on the left hand side and the function ${\bf B}_{u}$  are evaluated at 
$(y_{1}, v(y_{1}, z), z).$ After integrating over $(y_1,z)\in U\cap \{z\in L_{i} \cap G\}$ and summing over the finitely many coordinate charts employing a suitable partition of unity 
subordinate to the coordinate charts, we obtain 
\begin{equation*}
\int_{L_i\cap G}d{z}\int_{\ell_{z}}|\k_{z}|ds\leq \int_{M_i\cap T_2^{-1}(G)}{\mathbf B}_{u}\,
d\H^{n-1}.
\end{equation*}
By letting $i\->\00$ in this and using the Cauchy-Schwarz inequality on the right hand side, we deduce the desired estimate. 
\hfill$\Box$\\

We now proceed to derive the contradiction necessary to prove Proposition~\ref{mainprop}. First note that by Lemmas~\ref{lemma1} and~\ref{lemma2} we have that
\begin{equation}
\lim_{i\->\00}\L^{n-2}(B_1^{n-2}\setminus(D_{\e_i}\cap Q_{\e_i}))=0.
\label{van}
\end{equation}

For the rest of the proof let $T_2:B_{1}^{2} \times B_{1}^{n-2}\-> B_{1}^{n-2}$ be the orthogonal projection. By the defining property of $D_{\e_i}$ we have that
\begin{equation*}
\int_{\Ti \cap\{|u_{\e_i}|\leq \frac12\}}{\mathbf B}_{\e_i}^2|\nabla u_{\e_i}|\, d{x}\leq
\frac{\e_i}{c_3}\int_{B_1^2\times B_1^{n-2}} {\mathbf B}_{\e_i}^2 |\nabla u_{\e_i}|^2\, d{x}
\end{equation*}
which by \eqref{secdec2} tends to 0 as $i\->\00$. By the co-area formula it
then follows that
\begin{equation*}
\lim_{i\->\00}\int_{-\frac12}^{\frac12}dt\int_{\Ti\cap\{u_{\e_i}=t\}}
{\mathbf B}_{\e_i}^2\, d\H^{n-1}=0.
\end{equation*}

Now choose a generic $t \in (-1/2, 1/2)$ such that $\{u_{\e_i}=t\}$ is a $C^3$ hypersurface of
$B_1^2\times B_1^{n-2}$ for each $i=1, 2, 3, \ldots,$ 
\begin{equation*}
\lim_{i\->\00}\int_{\Ti\cap\{u_{\e_i}=t\}}{\mathbf B}_{\e_i}^2\, d\H^{n-1}=0
\end{equation*}
and 
\begin{equation*}
\liminf_{i\->\00} \, \H^{n-1}(\Ti\cap\{u_{\e_i}=t\})<\00.
\end{equation*}
This last requirement can be met since by the co-area formula and Fatou's lemma, 
\begin{eqnarray*}
&&\int_{-1/2}^{1/2} \liminf_{i \to \infty} \, {\mathcal H}^{n-1}(\Ti\cap\{u_{\e_i}=t\}) \, dt\\
&&\hspace{1.5in} \leq \limsup_{i\->\00}\int_{\Ti\cap\{|u_{\e_i}|\leq \frac12\}}|\nabla u_{\e_i}|\, d{x}\\
&&\hspace{1.5in} \leq \limsup_{i\->\00}\frac{\e_i}{c_3}\int_{B_2}|\nabla u_{\e_i}|^2\, d{x}<\00,
\end{eqnarray*}
Applying Lemma \ref{lemma3} with $u = u_{\e_{i}}$ and $G=D_{\e_i}\cap Q_{\e_i}$, we see in view of \eqref{van} 
that after passing to a subsequence without changing notation, there is a point ${z}_i\in D_{\e_i}\cap Q_{\e_i},$  
$i  =1, 2, 3, \ldots,$ such that 
\begin{equation}\label{flat}
\lim_{i\->\00}\int_{\ell^{i}_{{z}_i}}|\k^{i}_{{z}_i}|\, ds=0, 
\end{equation}
where $\ell^{i}_{z} = T_{2}^{-1}(z) \cap \{u_{\e_{i}} = t\}$ and $\k^{i}_{z}$ is the curvature of the curve $\ell^{i}_{z}$. 
Note that ${\ell}^{i}_{{z}_{i}}$ is the union of disjoint, connected, embedded planar
curves having at least one point in each of the disjoint balls $B^{2}_{\delta}({\mathbf p}_{j}) \times \{z_{i}\}$, $j=1, 2, \ldots, N,$  and no boundary point in $B_{1}^{2} \times B_{1}^{n-2}.$ Since $N \geq 3$, there must exist at least one connected component $\gamma_{i}$ of $\ell^{i}_{{ z}_{i}}$ such that as one moves along $\gamma_{i} \cap \overline{B_{1}^{2} \times B_{1}^{n-2}}$ from one end point  
to the other, or once around $\gamma_{i}$ in case $\gamma_{i}$ is closed, (a continuous choice of) the unit normal $\nu_{i}$ to $\gamma_{i}$ changes by a fixed amount depending only on the cone ${\mathbf C}$ and independent of $i.$ By integrating the derivative $\nu_{i}^{\prime}(s)$ along $\gamma_{i}(s),$ where $s$ is the arc length parameter, 
and using the fact that $\nu_{i}^{\prime}(s) = \k^{i}_{z_{i}}(s)$, we then obtain a fixed positive lower bound for $\int_{\gamma_{i}} |\kappa_{\gamma_{i}}(s)| ds$ independent of $i$, contradicting \eqref{flat}.
This completes the proof of Proposition~\ref{mainprop}.
\hfill$\Box$

\end{document}